\documentclass{article}
\usepackage[utf8]{inputenc}

\usepackage[english]{babel}

\usepackage[utf8]{inputenc}
\usepackage[T1]{fontenc}

\usepackage{amssymb}
\usepackage{stmaryrd}
\usepackage{amsmath,amsfonts}
\usepackage[tikz]{bclogo}
\usepackage[top=4.34cm]{geometry}

\usepackage{lmodern}
\usepackage{textcomp}
\usepackage{mathtools, bm}
\usepackage{amssymb, bm}
\usepackage[unq]{unique}
\usepackage{enumerate}
\usepackage{comment}
\usepackage[breaklinks]{hyperref}

\usepackage[numbers]{natbib}  

\usepackage[breaklinks]{hyperref}
\usepackage[numbers]{natbib}

\usepackage{amsthm}

\usepackage{cleveref}

\newtheorem{theorem}{Theorem}[section]
\newtheorem{lemma}[theorem]{Lemma}

\newtheorem{conjecture}[theorem]{Conjecture}

\newtheorem{Claim}[theorem]{Claim}

\newtheorem{Corollary}[theorem]{Corollary}

\theoremstyle{definition}

\newcommand{\cC}{\mathcal{C}}

\title{Progress towards the 1/2-Conjecture for the domination game}
\author{Julien Portier\footnote{\href{mailto:jp899@cam.ac.uk}{jp899@cam.ac.uk}, Department of Pure Mathematics and Mathematical Statistics (DPMMS), University of Cambridge, Wilberforce Road, Cambridge, CB3 0WA, United Kingdom} \and Leo Versteegen \footnote{\href{mailto:lvv23@cam.ac.uk}{lvv23@cam.ac.uk}, Department of Pure Mathematics and Mathematical Statistics (DPMMS), University of Cambridge, Wilberforce Road, Cambridge, CB3 0WA, United Kingdom}}
\date{}

\begin{document}

\maketitle

\begin{abstract}
The domination game is played on a graph $G$ by two players, Dominator and Staller, who alternate in selecting vertices until each vertex in the graph $G$ is contained in the closed neighbourhood of the set of selected vertices. Dominator's aim is to reach this state in as few moves as possible, whereas Staller wants the game to last as long as possible. In this paper, we prove that if $G$ has $n$ vertices and minimum degree at least 2, then Dominator has a strategy to finish the domination game on $G$ within $10n/17+1/17$ moves, thus making progress towards a conjecture by Bujt{\'a}s, Iršič and Klavžar.
\end{abstract}

\section{Introduction}

The domination game is a game for two players introduced by Bre{\v{s}}ar, Klav{\v{z}}ar and Rall in \cite{brevsar2010domination} and is played on a graph $G$ without isolated vertices. In this game, the two players, Dominator (D) and Staller (S), alternately select vertices of $G$. A vertex is considered to be \emph{dominated}, if itself or one of its neighbors has been selected by either player. Both players must always select a vertex that increases the number of dominated vertices. The game ends once all vertices of $G$ are dominated. Within the bound of these rules, Dominator wishes to minimize the number of moves played during the game, while Staller wants to achieve the opposite. We define \emph{the game domination number} $\gamma_{g}(G)$ of $G$ to be the number of moves the game takes to finish if Dominator moves first and both players play optimally. 

Let $n$ be the number of vertices of $G$. Of course, one should expect that the game lasts longer the larger $n$ is, and so it is a natural objective for research to bound $\gamma_g$ in terms of $n$. Specifically, Kinnersley, West and Zamani conjectured \cite{kinnersley2013extremal} that $\gamma_g(G)\leq 3n/5$ for all graphs without isolated vertices, and the second author has recently proved their conjecture \cite{versteegen2022domination}.

\begin{theorem}\label{thm:3-5}
    For any graph $G$ on $n$ vertices without isolated vertices, $\gamma_g(G)\leq 3n/5$.
\end{theorem}

 The bound $3n/5$ in \Cref{thm:3-5} is best possible, as one may see by considering the path on five vertices. What is more, if we think of the path $P_5$ as a single vertex to which we have appended two paths of length two, we can generalise this example by taking any graph $H$ and appending two paths of length two to each of its vertices. The resulting graph $G$ has $5\vert V(H)\vert$ vertices and satisfies $\gamma_g(G)=3\vert V(H)\vert$. One can generalise this further (see \cite{henning_loewenstein2017, versteegen2022domination}), but one property that all the resulting extremal graphs share is that they have many leaves, i.e., vertices of degree one. There is an extremal graph without a leaf, namely the cycle on five vertices, but in many ways, the five-cycle is a degenerate example that only works due to its small size, and unlike the $P_5$, it does not lend itself to generalisations. With this caveat in mind, it is natural to ask whether \Cref{thm:3-5} can be improved for the class of graphs without leaves. In \cite{bujtas20221}, Bujt{\'a}s, Iršič and Klavžar made the following two conjectures in this direction.

\begin{conjecture}\label{weak12Conjecture}
There exists a constant $c < 3/5$ such that every graph $G$ on $n \geq 6$ vertices and minimum degree at least $2$ satisfies $\gamma_g(G) \leq cn$.
\end{conjecture}

\begin{conjecture}\label{12Conjecture}
Let $G$ be a graph with $\delta(G) \geq 2$. Then $\gamma_g(G) \leq \lceil \frac{1}{2}n \rceil$.
\end{conjecture}

Of course, \Cref{12Conjecture} is a straightforward strengthening of \Cref{weak12Conjecture}, as it essentially sets the undetermined constant $c$ from \Cref{weak12Conjecture} to $1/2$, which is natural choice for two reasons. Firstly, it is simply the least constant for which there are not any counterexamples known already. For example, for any $n$ that does not have residue $3 \mod 4$, it is known \cite{kosmrlj2017Cycles} that the cycle $C_n$ satisfies $\gamma_g(C_n)=\lceil \frac{1}{2}n \rceil$. On the other hand, if we follow the heuristic notion that $\gamma_g$ should generally decrease if we add more edges to the graph under consideration, then we should also believe that cycles pose the worst case when proving an upper bound on $\gamma_g$ for graphs with minimum degree 2. 

In 2016, long before \Cref{thm:3-5} was established in full, Henning and Kinnersley proved \cite{henning2016domination} that \Cref{thm:3-5} holds for graphs of minimum degree 2, which is the reason that \Cref{weak12Conjecture} asks for a constant $c< 3/5$. To the best of our knowledge, no improvement on the constant has been made since, whence \Cref{weak12Conjecture} is still open. In this paper, we prove the following.

\begin{theorem}\label{MainThm}
For every graph $G$ with minimum degree at least $2$, we have $\gamma_g(G) \leq \frac{10}{17}n+\frac{1}{17}$.
\end{theorem}

 Since $10n/17+1/17<3n/5$ for all $n>5$, \Cref{MainThm} confirms \Cref{weak12Conjecture}.

\section{The proof}

\subsection{The main idea of the proof}

To prove \Cref{MainThm}, we employ the \emph{discharging method} introduced by Bujt{\'a}s in \cite{discharging_method}. The core concept of this method is that at any given time during the game, every vertex has a \emph{color}, which has in turn a certain number of \emph{points} associated with it. As the game progresses, the sum $\pi$ of all points associated to vertices decreases until it reaches zero. If we can prove that the rate of decrease per move is sufficiently large, then that will allows us to bound the total number of moves in terms of the total number of points at the start of the game, which in turn only depends on the number of vertices $n$ in the graph $G$ under consideration.

Concretely, we use the colours \emph{white, blue, orange} and \emph{red}. Let $t$ be the number of moves that have been played so far and let $A$ be the set of vertices that have been selected by either player up to and including move $t$. The colour of a vertex is determined by the following rules.

\begin{itemize}
    \item A vertex is white if and only if it is not dominated by $A$ (i.e, is not itself contained in $A$ and not adjacent to a vertex in $A$).
    \item A vertex is red if and only if it and all its neighbours are dominated by $A$.
    \item Let $W$ be the set of white vertices as determined by $A$. An isolated vertex in the induced subgraph $G[W]$ is called a \emph{single}, while a component of exactly two vertices in $G[W]$ is called a \emph{pair}. An individual member of a pair is called a \emph{double}. A vertex is orange if and only if it is dominated by $A$, has exactly one white neighbour, and this neighbour is a single or a double.
    \item All other vertices are blue.
\end{itemize}

The number of points $\pi_v(t)$ associated with a vertex $v$ after move $t$, depends only on the color of $v$ after move $t$. Namely, $\pi_v(t)$ is 20, 10, 7, or 0 if $v$ is white, blue, orange, or red, respectively. With this we define the total number of points after move $t$ as $\pi(t)=\sum_{v\in V(G)} \pi_v(t)$.


At the start of the game, all vertices are white, which means that $\pi(0)=20n$, and the game ends after move $T$ if all vertices are red, which is to say that $\pi(T)=0$. Therefore, if D can play in a way such that we have $\pi(T)\leq 20n-34T$, we would know that $T\leq 10n/17$. One can show that this is the case as long as $n>5$, but doing so requires a disproportional amount of additional effort relative to the only very slightly worse bound $\pi(T)\leq 10n/17+1/17$. To show the latter bound, we divide the game for the purpose of Dominator's strategy into two phases. The first phase lasts as long as D can play in a way that guarantees the game to arrive at a point where $t$ moves have been played, $\pi(t)\leq \pi(0)-34t$, and D is to move next. 

Once this phase is over, one can conclude that all the remaining white vertices lie in very specific configurations. The problem with these configurations is that they lead to a form of zugzwang, where neither player would like to make the first move. A good example for such a configuration is a cycle on five white vertices. Indeed, it is an easy exercise to show that if the domination game is played on a $C_5$ and D makes the first move, then the game will take three moves, whereas it will only take two moves if S goes first. Happily, if the five-cycle is part of a larger graph $G$, it is S who is in zugzwang after those three moves have been played, and there is no way for them to get out of it for the rest of the game. Thus, the idea is that D can afford to reduce $\pi$ by less than 34 on average for a few moves in exchange laying the burden of the initiative on S.

\subsection{The active phase}

Before we reveal the strategy that D should play, we first make a general observation about the point system.

\begin{Claim}\label{claim:min-points}
    If the game is not over after $t$ moves, then $\pi(t)-\pi(t+1)\geq 20$.
\end{Claim}

\begin{proof}
Suppose that the vertex $v$ is played on move $t+1$. If $v$ was blue before move $t+1$, then $v$ has a neighbour $u$ which was white then. After move $t+1$ on the other hand, $v$ is red and $u$ is not white anymore, so that $\pi(t)-\pi(t+1)\geq 10+(20-10)= 20$. 
If $v$ was coloured orange, then $v$ has a neighbour $u$ which was a double or a single before move $t+1$. In particular, $u$ was white and has at most one white neighbour, whence it will be orange or red after move $t+1$. It follows that $\pi(t)-\pi(t+1)\geq 7+(20-7)= 20$. 
Finally, if $v$ was white itself, $\pi(t)-\pi(t+1)\geq 20$ because $v$ is red after move $t+1$.
\end{proof}

We say that the game is in a \emph{stable state} after move $t$, if $\pi(0)-\pi(t)\geq 34t$, and at least one of the following holds.
\begin{itemize}
    \item Either, $t$ is even, i.e., D is to move next, or
    \item there are no white vertices left, i.e., the game is over.
\end{itemize}
Whenever the game is not over but in a stable state, it must be Dominator's turn. Suppose that the game is in a stable state after move $t$. If D has a strategy that allows him to advance the game, possibly over the course of several moves, to another stable state after the current one, D plays according to this strategy. Letting $t_1$ denote the smallest integer so that D does not have such a strategy after move $t_1$, we can make a series of powerful observations about the state of the game after move $t_1$. 

For example, we know that there can be no vertex $v$ such that $\pi(t_1)-\pi(t_1+1)\geq 48$ if D plays $v$ on move $t_1+1$. Indeed, if the game is over after D plays $v$, the game has reached another stable state after move $t_1+1$. If it is not, \Cref{claim:min-points} tells us that any subsequent move by S will reduce $\pi$ further by at least 20 so that $\pi(t_1)-\pi(t_1+2)\geq 68=2\cdot 34$, which means the game will reach a stable state since $t_1+2$ is even. Likewise, there can never be a strategy for D to play that guarantees that $\pi(t_1)-\pi(t_1+2)\geq 68$ or $\pi(t_1)-\pi(t_1+3)\geq 116$, even if it begins with a move that decreases $\pi$ by less than 48. In the proofs of all of the following claims, the assumption that the claim does not hold can be used to provide D with a strategy that leads to one of these three outcomes and thus yields a contradiction. 

\begin{Claim}\label{claim:no-triple}
    No white vertex has three white neighbours and no blue vertex has four white neighbours.
\end{Claim}

\begin{proof}
    Suppose towards a contradiction that there is a white vertex $x$ with at least three white neighbours after $t_1$ moves. D may play this vertex, and since all neighbours of $x$ are now dominated, each of them has at least 10 points less associated with them than before. Furthermore, $x$ is now red and has 0 point associated with it so that $\pi(t_1)-\pi(t_1+1)\geq 50$. Similarly if a vertex $v$ is blue and has four white neighbours, D can also reduce $\pi$ by at least 50 by playing $v$.
\end{proof}

Let $H^*$ be the subgraph of $G$ in which we remove every connected component that is a four- or five-cycle and has only white vertices. Let further $H$ be the subgraph of $H^*$ that is induced by the set of white vertices after $t_1$ moves. By \Cref{claim:no-triple}, $H$ has no vertex of degree three or higher, whence $H$ must be a union of disjoint paths and cycles.

\begin{Claim}
    $H$ contains no path $xyz$ such that $x$ has no white neighbour other than $y$ and $z$.
\end{Claim}

\begin{proof}
    If such a path existed, then D could simply play $y$ to decrease $\pi$ by 50.
\end{proof}

\begin{Claim}
    $H$ contains no cycle on nine or more vertices.
\end{Claim}

\begin{proof}
Suppose that $H$ has a $k$-cycle $x_1x_2\ldots x_kx_1$ for some $k \geq 9$. In this case D may play $x_5$, thus decreasing $\pi$ by 46. If S replies by playing a vertex $z$ that is adjacent to both a vertex in $\{x_1,x_2,x_3\}$ and $\{x_7,x_8,x_9\}$, then $\pi$ decreases further by at least 30, so we may assume without loss of generality that $y$ has no neighbour in $\{x_1,x_2,x_3\}$. In this case, D can play $x_2$ next, which decreases $\pi$ by 60 since both $x_3$ and $x_4$ become red. By \Cref{claim:min-points}, we have $\pi(t_1)-\pi(t_1+3)\geq 126$.
\end{proof}

\begin{Claim}
    $H$ contains neither a eight-cycle nor a seven-cycle.
\end{Claim}

\begin{proof}
Suppose that $H$ has an eight-cycle $x_1x_2\ldots x_8x_1$. D can play $x_1$, which decreases $\pi$ by 40, and if the vertex $z$ that S plays in response is adjacent to two white vertices, then $\pi$ decreases by at least 30. We may therefore assume that $z$ has only one white neighbour before Staller's move and thus without loss of generality that $z$ is not adjacent to either of $x_3$ and $x_4$. If $z$ is not adjacent to $x_5$ either, then D can play $x_4$ on move $t_1+3$, which decreases $\pi$ by another 60 points, whence we would have $\pi(t_1)-\pi(t_1+3)\geq 120$. On the other hand, if $z$ is adjacent to $x_5$, then both $x_5$ and $x_8$ will be orange after D plays $x_4$ on move $t+3$, whence $\pi(t_1)-\pi(t_1+3)\geq 116$.

The proof that $H$ cannot have a seven-cycle is very similar and we leave it as an exercise.
\end{proof}

\begin{Claim}
    $H$ contains no six-cycle.
\end{Claim}
\begin{proof}
    Suppose that $H$ contains a six-cycle $x_1\ldots x_6x_1$. D can play $x_1$, and if after Staller's move any of the vertices $x_3,x_4$ and $x_5$ is still white, D can play $x_4$ next and we have $\pi(t_1)-\pi(t_1+3)\geq 120$. Of course, if after move $t_1+2$, none of $x_3,x_4$ and $x_5$ are white anymore, we even have $\pi(t_1)-\pi(t_1+2)\geq 120$.
\end{proof}

\begin{Claim}
    $H$ contains no four-cycle.
\end{Claim}
\begin{proof}
    Suppose that $H$ contains a four-cycle $x_1\ldots x_4x_1$. Because $H$ is a subgraph of $H^*$, one of the four vertices, $x_1$ say, must have another neighbour $u$, which must be blue. D can play $x_3$, which makes $x_3$ red as well as $x_2$ and $x_4$ orange. We distinguish which colour $x_1$ and $u$ have after S played some vertex $v$ on move $t_1+2$. If $x_1$ is red, then $x_2$ and $x_4$ are red as well, and we must have $\pi(t_1)-\pi(t_1+2)\geq 80$. If $x_1$ is not red, then it must be white, and thus $u$ cannot be red either. In particular, $u\neq v$. If $u$ is orange, we can use the fact $v$ is now red, and that $v$ must have been white or had a white neighbour before move $t_1+2$ to see that $\pi(t_1)-\pi(t_1+2)\geq 69$.

    This leaves the possibility that $u$ is still blue after move $t_1+2$. In this case, $u$ must have another white neighbour $y$ besides $x_1$, and D can play $u$ after which $u, x_1,x_2$ and $x_4$ must be red and $y$ is not white anymore so that $\pi(t_1)-\pi(t_1+3)\geq 120$ by \Cref{claim:min-points}.
\end{proof}

We know now that $H$ must consist only of isolated vertices, isolated edges and five-cycles, but the structure of $H$ can be even more rigidly determined. Indeed, we will show that $H$ does in fact not contain any five-cycles or isolated vertices and that the edges of $H$ lie in highly specific configurations when viewed in $G$. In the following proofs, it is understood that whenever we play a vertex that is adjacent to a double, the double becomes orange or red but never blue.

\begin{Claim}\label{claim:no-two-c5}
    If two vertices $x_1,y_1\in H$ are part of (not necessarily distinct) five-cycles, they cannot have a common blue neighbour.
\end{Claim}
\begin{proof}
    Suppose that there are two five-cycles $x_1\ldots x_5x_1$ and $y_1\ldots y_5y_1$ in $H$ and a vertex $u$ that is adjacent to both $x_1$ and $y_1$. If $y_1=x_i$ for some $i\neq 1$, there must be some $j\in [5]$ such that $\{x_1,x_i\}\subset \{x_j\}\cup N(x_j)$. Therefore, by \Cref{claim:no-triple}, $u$ will become red or orange after D plays $x_j$ so that $\pi$ decreases by at least 49. Let us therefore assume that the two cycles are disjoint. 
    
    Suppose now that $u$ has a third white neighbour $z$, noting that $u$ cannot have a fourth white neighbour by \Cref{claim:no-triple}. D can play $u$ which reduces $\pi$ by at least 40. If S replies by playing a vertex that has a neighbour in $\{x_2,x_3,x_4,x_5\}$ and $\{y_2,y_3,y_4,y_5\}$, then $\pi(t_1)-\pi(t_1+2)\geq 70$ and we are done. Otherwise, we can assume without loss of generality that $\{x_2,x_3,x_4,x_5\}$ are all still white after move $t_1+2$, and D can play $x_3$ which colours $x_2$ and $x_3$ red as well as $x_1$ and $x_4$ orange, thus reducing $\pi$ by at least 56 so that $\pi(t_1)-\pi(t_1+3)\geq 116$.

    We are left with the possibility that $x_1$ and $y_1$ are the only white neighbours of $u$. In this case, D may play $x_1$, which reduces $\pi$ by 46. If S replies by playing a vertex that is adjacent to at least two of the vertices $\{x_3,x_4,y_1,y_2,y_5\}$, they reduce $\pi$ further by at least 30. Should S play a vertex that does not dominate either of $x_3$ or $x_4$, D can play $x_3$ next so that $\pi(t_1)-\pi(t_1+3)\geq 120$. Lastly, if S plays a vertex that does not dominate either of $y_1,y_2$ and $y_5$, D can play $y_1$ next so that $\pi(t_1)-\pi(t_1+3)\geq 122$.
\end{proof}

\begin{Claim}\label{claim:blue-deg2}
    No (blue) vertex has three white neighbours.
\end{Claim}

\begin{proof}
    Let $u$ be a blue vertex with three white neighbours $x_1,y$ and $z$. If all three of them are singles or doubles, D can play $u$ to decrease $\pi$ by at least 49. Let us therefore assume that $x_1$ is part of a five-cycle $x_1\ldots x_5x_1$ in $H$. By \Cref{claim:no-two-c5}, $y$ and $z$ must be singles or doubles. If $y$ is a single, D can still play $u$ to reduce $\pi$ by at least 53 so let us assume that $y$ is adjacent to another double $y'$, which must have itself a second neighbour, which we denote by $v$. How D should proceed depends on the number and type of the white neighbours of $v$.

    Suppose first that $v$ has only one white neighbour, D can play $y$ to reduce $\pi$ by 47. If S replies by playing a vertex that dominates $x_1$ or at least two white vertices overall, $\pi$ decreases further by at least 23 and we are done. Otherwise, D can play one of $x_1,x_2$ or $x_5$ next to reduce $\pi$ by at least 49.

    Suppose next that $v$ has exactly two white neighbours and that one of them (the one that is not $y'$), which we call $x_1'$, is part of a five-cycle $x_1'\ldots x_5'x_1'$ of white vertices. In this scenario, D can reduce $\pi$ by 49 by playing $x_1'$.

    If neither of the above assumptions is true, we know by \Cref{claim:no-two-c5} that $v$ is adjacent to a vertex $z'\neq y'$ which is either a single or a double. If $z'=y$, D can play $v$ to reduce $\pi$ by 50 right away. Else, D can play $u$ on move $t_1+1$ to reduce $\pi$ by at least 46. As usual, we may assume that the vertex that S plays in response, dominates only one white vertex. We can also assume that this white vertex is either $y'$ or $z'$ as otherwise D could subsequently reduce $\pi$ by at least 50 by playing $v$. In particular, $x_2,x_3$ and $x_4$ are still white after move $t_1+2$, which means that D can play $x_3$ on move $t_1+3$, thus reducing $\pi$ by at least 56.
\end{proof}

\begin{Claim}
    $H$ contains no five-cycle.
\end{Claim}

\begin{proof}
Suppose that $H$ has a five-cycle $x_1\ldots x_5x_1$. We know by definition of $H^*$ and $H$ that one of the vertices in our cycle, $x_1$ say, must have another neighbour $u$, which must be blue. By \Cref{claim:no-triple} and \Cref{claim:blue-deg2}, $u$ can have zero or one further white neighbours. In either case, D can play $x_1$, which reduces $\pi$ by at least 49.
\end{proof}

\begin{Claim}\label{claim:no-double-orange}
    No double is adjacent to an orange vertex.
\end{Claim}
\begin{proof}
    Suppose $xy$ is a pair and that $x$ is adjacent to an orange vertex $u$. Since no vertex has degree 1, $y$ must have another neighbour $v$, which cannot be white and must be distinct from $u$. Furthermore, $v$ cannot have three white neighbours by \Cref{claim:blue-deg2} and if D plays $x$, then $\pi$ will decrease by at least 50.
\end{proof}

\begin{Claim}\label{claim:no-single-blue-orange}
    No single is adjacent to both a blue and an orange vertex.
\end{Claim}
\begin{proof}
    By playing the blue neighbour of the single, D would decrease $\pi$ by at least 50.
\end{proof}

\begin{Claim}\label{claim:no-blue-pair}
    No blue vertex is adjacent to both doubles of a pair.
\end{Claim}
\begin{proof}
    By playing such a blue vertex, D would decrease $\pi$ by at least 50.
\end{proof}

\begin{Claim}\label{claim:no-double-two-blues}
    No double is adjacent to two blue vertices.
\end{Claim}
\begin{proof}
    By playing the double, D would decrease $\pi$ by at least 49.
\end{proof}

\begin{Claim}\label{claim:no-two-singles}
    No (blue) vertex is adjacent to two singles.
\end{Claim}
\begin{proof}
    By playing such a blue vertex, D would decrease $\pi$ by at least 56.
\end{proof}

\begin{Claim}\label{claim:no-single-blue}
     No single has a blue neighbour.
\end{Claim}

\begin{proof}
    Suppose towards a contradiction that there exists a single $x$ with a blue neighbour $u$. Since the degree of $x$ is at least 2, $x$ must have another neighbour $v$, and by \Cref{claim:no-single-blue-orange}, $v$ must be blue. By \Cref{claim:no-two-singles}, $u$ and $v$ must both be adjacent to a double, which we call $y$ and $z$ respectively. 
    
    Assume first that $y$ and $z$ are not adjacent. In this case, D can first play $y$ which decreases $\pi$ by 46. If S replies by playing a vertex that is adjacent to $x$ or $z$, that will decrease $\pi$ by at least 23 and we are done. In any other case, D may play $v$ next to decrease $\pi$ by at least 50, so that $\pi(t_1)-\pi(t_1+3)\geq 116$.

    If $y$ and $z$ are adjacent, D can still decrease $\pi$ by 46 by playing $y$. By \Cref{claim:no-double-two-blues} and \Cref{claim:no-double-orange}, $y$ and $z$ have no other neighbours than each other and $u$ and $v$ respectively. Therefore, \Cref{claim:no-double-orange} is still valid after this move by D. Hence, if S plays an orange vertex, then $\pi$ will decrease by at least 34 due to \Cref{claim:no-single-blue-orange}. If S plays a single, all of its orange neighbours become red and all of its blue neighbours become orange so that $\pi$ reduces by at least 26. If S plays a blue vertex, $\pi$ decreases by at least 36 and if S plays a double, $\pi$ decreases by 46. In any of these cases, we have $\pi(t_1)-\pi(t_1+2)\geq 72$.
\end{proof}

One can conclude from \Cref{claim:no-double-orange} and \Cref{claim:no-single-blue} that S would have no good move if it was their turn rather than Dominator's. Often, this will allow D to bring S into a zugzwang of sorts. The next claim formalizes this.

\begin{Claim}\label{claim:zugzwang}
    Suppose that after move $t\geq t_1$, no double is adjacent to an orange vertex and no single is adjacent to both a blue and an orange vertex. Then move $t+1$ will decrease $\pi$ by at least 34.
\end{Claim}
\begin{proof}
We denote the vertex played on move $t+1$ by $x$. If $x$ is a double, $\pi$ decreases by 46, and if $x$ is blue, $\pi$ decreases by 36. If $x$ is a single, it must have at least two orange neighbours by \Cref{claim:no-single-blue} and $\pi$ decreases by at least 34. The same holds if $x$ is orange itself.
\end{proof}

As a simple consequence of \Cref{claim:zugzwang}, we obtain the following.

\begin{Claim}\label{claim:no-single}
    No orange vertex and no single exists.
\end{Claim}
\begin{proof}
    If an orange vertex existed, D could reduce $\pi$ by 34 by playing this vertex without introducing a new orange or blue vertex, whence \Cref{claim:zugzwang} would be applicable with $t=t_1+1$ (recall \Cref{claim:no-double-orange}).

    Since no orange vertex exists, no single can exist either by \Cref{claim:no-single-blue}.
\end{proof}

Another consequence of \Cref{claim:zugzwang} is that the isolated four-cycles we removed in our definition of $H^*$ do not actually pose a problem.

\begin{Claim}\label{claim:no-iso-c4}
After move $t_1$, $G$ has no isolated four-cycle with only white vertices.
\end{Claim}
\begin{proof}
If $G$ had such a cycle, D could just play any vertex to reduce $\pi$ by 46, after which \Cref{claim:zugzwang} is applicable.
\end{proof}

\begin{Claim}\label{claim:double-c6}
    Every pair $xy$ must lie inside a six-cycle $uxyvy'x'u$ such that $x'$ and $y'$ are white. Furthermore, the six-cycle is unique up to relabelling of its constituent vertices.
\end{Claim}

\begin{proof}
To see that any six-cycle as in the claim must be unique, note that by \Cref{claim:blue-deg2}, $u$ has no other white neighbours than $x$ and $x'$, and $v$ has no white neighbours other than $y$ and $y'$. 

Let now $xy$ be a pair. By \Cref{claim:no-double-orange} and \Cref{claim:no-double-two-blues}, $x$ and $y$ have both exactly one blue neighbour, which we denote by $u$ and $v$ respectively. These in turn must have exactly one more white neighbour each, which we call $x'$ and $y'$ respectively. By \Cref{claim:no-single}, both $x'$ and $y'$ are doubles, by \Cref{claim:no-blue-pair}, $y'\neq x$ and $x'\neq y$ and by \Cref{claim:no-double-two-blues}, $x'\neq y'$.

If $x'$ and $y'$ were adjacent, there would be nothing left to show, so we may suppose towards a contradiction that they are not. In this situation, D may play $x$ to decrease $\pi$ by 46. We can argue as in the proof of \Cref{claim:no-single} that if S plays any vertex other than $u$ or $v$, $\pi$ will decrease by at least 34. If S does play one of $u$ and $v$, without loss of generality the former, then $\pi$ decreases by 20 and D can play $y'$ next, thus decreasing $\pi$ by at least 50 so that we have $\pi(t_1)-\pi(t_1+3)\geq 116$.

\end{proof}

\subsection{The reactive phase}

We now know that after move $t_1$, all remaining white vertices in $G$ have degree 2 and lie either in an isolated five-cycle of white vertices or in a six-cycle as described in \Cref{claim:double-c6}. Let $k$ be the number of configurations of either type. The following claim asserts that D can rapidly end the game from here on out.

\begin{Claim}\label{claim:reactive-phase}
    D can play in a way such that the game ends after at most $t_1+2k+1$ moves.
\end{Claim}
\begin{proof}
    D can ensure that S plays at most one vertex in any of the $k$ cycles, from which the claim follows. Indeed, whenever S plays a vertex $x$ in one of the remaining cycles $C$, and there are white vertices left in $C$, D can make all of them red by playing a vertex that is opposite to $x$ in $C$. If there are no white vertices left in $C$ after Staller's move as well as on move $t_1+1$, D can simply play an arbitrary vertex.
\end{proof}

\subsection{Putting everything together}

\begin{proof}[Proof of \Cref{MainThm}.]
We have $\pi(0)=20n$, and by definition of $t_1$, $\pi(t_1)\leq \pi(0)-34t_1$. Furthermore, since each of the $k$ remaining five- and six-cycles carries 100 points, we know that $\pi(t_1)=100k$. Therefore, we obtain $100k\leq 20n-34t_1$ or 

\begin{align}\label{eq:t1}
    t_1\leq \frac{10}{17}n-\frac{50}{17}k.
\end{align}

If $k=0$, then the game is over after move $t_1$, and there is nothing left to show. Otherwise, we apply \Cref{claim:reactive-phase} to see that D can conclude the game after $t_1+2k+1$ moves, which by \eqref{eq:t1} is at most $10n/17+1/17$.
\end{proof}

\section{Concluding remarks}
 Before we discuss whether the proof of \Cref{MainThm} is amenable to further improvement, we point out that instead of letting D reduce $\pi$ greedily (as fast as possible), one can extract a more explicitly formulated strategy for D from the Claims \ref{claim:no-triple} through \ref{claim:reactive-phase}. Indeed, whenever it is Dominator's turn, one looks for the first one among these claims that is not yet satisfied, and lets D play according to the instructions in the proof of that claim. If all of the claims are satisfied, then one can move to play according to the explicit strategy of the reactive phase.\\

\begin{figure}[htbp]\centering
    			\includegraphics{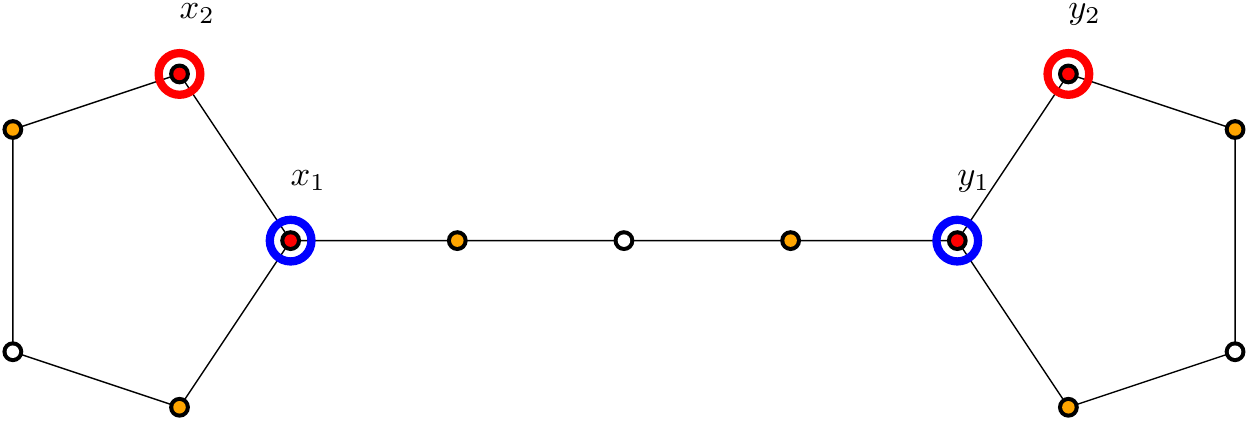}
    			\caption{The graph $G_{13}$. After D plays the vertices with encircled in blue, and S plays the vertices encircled in red, the vertices are coloured as shown here.}
    			\label{fig:G13}
\end{figure}

We believe that with a more careful case analysis the method used in this paper could yield a better constant $C < \frac{10}{17}$ in \Cref{MainThm}. However, it seems that a substantially different approach would be required to improve the bound beyond $\frac{7}{13}n$, as the example $G$ of multiple disjoint copies of the graph $G_{13}$ in \Cref{fig:G13} shows.

According to the strategy that the Claims \ref{claim:no-triple} to \ref{claim:reactive-phase} describe on $G$, D will start by playing the vertex $x_1$ or $y_1$ in one of the copies of $G_{13}$, to which S can reply by playing $x_2$ or $y_2$ in the same copy. After all of those four vertices have been played in every copy of $G_{13}$, three more vertices will have to be played in each of them. Overall, 7 out of 13 vertices will be played in every instance of $G_{13}$ over the course of the game.

We will now make some remarks on conjectures that are related to \Cref{12Conjecture}. Firstly, the proof of \Cref{MainThm} can be adapted to yield progress on the following conjecture by Rall, which was first mentioned in \cite{james2019domination}.

\begin{conjecture}\label{RallConjecture}
If $G$ has a Hamiltonian path, then $\gamma_g(G) \leq \lceil \frac{1}{2}n \rceil$.
\end{conjecture}

Indeed, we can show the following as an easy corollary to the proof of \Cref{MainThm}.

\begin{theorem}
If $G$ has a Hamiltonian path, then $\gamma_g(G) \leq \lceil \frac{10}{17}n \rceil$.
\end{theorem}

\begin{proof}[Sketch proof.]
The worst case is that $G$ has two leaves $l_1$ and $l_2$. In this situation, D can play the owner of $l_1$ in the first move to decrease $\pi$ by 50. If $l_2$ is not red before Dominator's next move, its unique neighbour cannot be red either, so playing it decreases $\pi$ by at least 27. Therefore, we have $\pi(0)-\pi(4)\geq 117$. If $l_2$ becomes red before the second move, we have $\pi(0)-\pi(2)\geq 60$. In either case, it is easy to check that if D plays as in the proof of \Cref{MainThm} from there on out, the game will last at most $10n/17+21/34\leq \lceil \frac{10}{17}n\rceil$ moves.
\end{proof}

As a relaxation of \Cref{12Conjecture}, Bujt{\'a}s, Iršič and Klavžar \cite{bujtas20221} posed the following conjecture.

\begin{conjecture}\label{12UniversalConstant}
There exists a universal constant $C$ such that for every graph $G$ with $\delta(G) \geq 2$, we have $\gamma_g(G) \leq \frac{1}{2}n +C$.
\end{conjecture}

It turns out however, that \Cref{12Conjecture} and \Cref{12UniversalConstant} are almost equivalent in an even stronger sense than it seems. To see this, we consider a variant of the domination game, which can be seen as a generalisation of most other variants that have been studied. The \emph{transversal game} was defined in \cite{bujtas2016transversal} as a game played on a hypergraph $H$, in which the two players, Edge-hitter and Staller, alternately select a vertex from $H$, with the rule that each newly selected vertex must \emph{hit}, i.e., be contained in, at least one edge that does not intersect the set of previously selected vertices. The game ends when the set of selected vertices becomes a transversal in $H$, i.e., when every edge of $H$ intersects the set of selected vertices. Edge-hitter aims to end the game in as few moves as possible, whereas Staller wants the game to last as long as possible. The \emph{game transversal number} $\tau_g(H)$ of $H$ is the number of moves played if Edge-hitter starts the game and both players play optimally. In \cite{portier2022proof}, the authors proved the following result.

\begin{lemma}\label{copylemmageneral}
Let $\cC$ be a set of hypergraphs closed under taking multiple disjoint copies of a hypergraph in $\cC$ and suppose that there exists $c>0$ such that for every hypergraph $H \in \cC$ on $n$ vertices, $\tau_{g}(H) \leq (c+o(1))n$. Then we also have $\tau_{g}(H) \leq cn+1$ for every $H\in \cC$ on $n$ vertices. 
\end{lemma}

The closed neighbourhood hypergraph of a graph $G$ is defined as the hypergraph $H_G$ with vertex set $V(H_G) = V(G)$ and hyperedge set $E(H_G) = \{ N_G[x] | x \in V(G) \}$ consisting of the closed neighbourhoods of vertices in $G$. As remarked in \cite{bujtas2016transversal}, the domination game played on a graph $G$ can be seen as a special instance of the transversal game being played on the closed neighbourhood hypergraph $H_G$ of $G$. By taking $\cC =\{H_G:G \text{ is a graph of minimum degree at least 2}\}$ and $c=1/2$, we obtain the following as a direct corollary of \Cref{copylemmageneral}.

\begin{Corollary}\label{copylemma12}
Suppose there exists a strategy for D to finish the domination game on each graph $G$ of minimum degree at least $2$ in at most $\frac{1}{2}(1+o(1))n$ moves. Then there also exists a strategy for D to finish the domination game on each graph $G$ of minimum degree at least $2$ in at most $\frac{1}{2}n+1$ moves.
\end{Corollary}

Thus, if \Cref{12UniversalConstant} is true, we can actually choose $C=1$.

\section*{Acknowledgement}

The authors would like to thank Béla Bollobás for his valuable comments.
\bibliographystyle{abbrvnat}  
\renewcommand{\bibname}{Bibliography}
\bibliography{bibliography}

\end{document}